\documentstyle[english,12pt,fleqn]{amsart}
\newfont{\para}{cmbsy10 scaled 1200}
\def\PS{{\para x}\ }

\newcommand{\EP}{{\mathbb  P}}

\newcommand{\EN}{{\mathbb  N}}

\newcommand{\EZ}{{\mathbb  Z}}
\newcommand{\EQ}{{\mathbb Q}}

\newcounter{aufgabecount}

\bigskipamount\medskipamount

\begin{document}
\parindent0em
\renewcommand{\labelenumi}{\alph{enumi})}

\title[Rational curves and ampleness properties of the 
tangent bundle] 
{Rational curves and ampleness properties of the 
tangent bundle of algebraic varieties}

\author{Fr\'ed\'eric Campana}
\address{Fr\'ed\'eric Campana\\
       Universit\'e de Nancy I   \\
       D\'epartement de Math\'ematiques \\
       B.P. 239 \\
       54506 Vandoeuvre les Nancy Cedex \\
       France}
\email{frederic.campana@@iecn.u-nancy.fr}

\author{Thomas Peternell}
\address{Thomas Peternell \\
       Mathematisches Institut \\
       Universit\"at Bayreuth \\
       D - 95440 Bayreuth \\
       Germany}
\email{peternel@@btm8x1.mat.uni-bayreuth.de}
\thanks{Parts of this paper were written up during a stay of one of the authors
at MSRI. He would like to thank the institute for support and the excellent
working conditions. Research at MSRI is supported on part by NSF grant
DMS-9022140.}

\begin{abstract}
The purpose of this paper is to translate positivity properties of the
tangent bundle (and the anti-canonical bundle) of an algebraic
manifold into existence and movability properties of rational curves
and to investigate the impact on the global geometry of the manifold
$X$. Among the results we prove are these:

\quad If $X$ is a projective manifold, and ${\cal E} \subset T_X$ is an
ample locally free sheaf with $n-2\ge rk {\cal E}\ge n$, then $X
\simeq \EP_n$.

\quad Let $X$ be a projective manifold.  If $X$ is rationally connected,
then there exists a free $T_X$-ample family of (rational) curves.  If
$X$ admits a free $T_X$-ample family of curves, then $X$ is rationally
generated.
\end{abstract}

\maketitle

\section*{Introduction}

The purpose of this paper is, vaguely speaking, to translate positivity
properties of the tangent bundle (and the anti-canonical bundle) of an algebraic
manifold into existence and movability properties of rational curves and to
investigate the
impact on the global geometry of the manifold $X$. This study falls into
two parts:
\begin{itemize}
\item[(1)] a biregular one in which ampleness of some subsheaf of the
tangent bundle $T_X$ is
assumed; one then expects a classification,
\item[(2)] a birational one, where the ampleness of $T_X$ is assumed only generically; one
then expects $X$ to be rationally connected.
\end{itemize}
In 1979 in his famous solution of the Hartshorne-Frankel conjecture,
S. Mori proved that a projective manifold $X$ whose tangent bundle $T_X$ is
ample, must be  the projective  space $\EP_n$. Previously, Mori and Sumihiro
showed that $\EP_n$ is the only manifold admitting a vector field  vanishing
along an ample divisor (and this is an Euler vector field if $n \ge 2)$.
This was generalised by J. Wahl to the extent that if
$T_X \otimes L^{-1}$ has a non-zero section with $L$ an ample line bundle (not
necessarily effective), then $X \simeq \EP_n$.

\bigskip

Although these two theorems look at first quite different in nature, they
might be only special cases of a more general theorem.

\bigskip

{\bf Question:} \hspace{0,3 cm} Let $X$ be a projective manifold,
${\cal E} $ an ample locally free sheaf,
${\cal E} \subset T_X.$ Is $X \simeq \EP_n$?

\bigskip

We shall prove

\bigskip
{\bf Theorem (1.1)} \hspace{0,3cm} {\it The question has a positive answer
if $rk {\cal E}$ is
$n-2$, $n-1$ or $n.$}

\bigskip

A manifold is said to be rationally connected if two general points can be
joined
by a chain of rational curves. Examples are Fano manifolds as shown by
Campana and Koll\'ar-Miyaoka-Mori. We ask whether it is possible to characterize
rational connectedness by a weak positivity property of the tangent bundle.
Philosophically speaking we require that $T_X$ should be ample on sufficiently
many curves; more precisely $T_X$ should be ample on the generic curve of a
{\it free}
family (2.5).

We shall prove in this spirit (for definitions see (2.1))

\bigskip
{\bf Theorem (2.7)} \hspace{0,3cm} {\it Let $X$ be a projective manifold.
\begin{itemize}
\item[(1)] If $X$ is rationally connected, then there exists a free
$T_X$-ample family
of (rational) curves.
\item[(2)] If $X$ admits a free $T_X$-ample family of curves, then $X$ is
rationally generated.
\end{itemize} }

\bigskip (1) is actually an easy consequence of [KoMiMo92]. A projective
variety $X$ is rationally
generated if every variety dominated by $X$ is uniruled, see sect.2.
This notion was introduced in [Ca95]); we expect that rationally generated
varieties are
actually rationally connected; this holds in dimension 3 by [KoMiMo92].

Important special cases of free family are given by complete intersections
of hyperplanes.
To make this precise we recall the following

\bigskip

{\bf Definition} (Miyaoka) \hspace{0,3 cm} Let $X$ be a $n-$dimensional
projective manifold, $H_1,..., H_{n-1}$ ample divisors on $X$. Let $ {\cal E}$
be a vector bundle on $X.$ Then $ {\cal E}$ is said to be generically
$(H_1,..., H_{n-1})$ -ample (nef) if for $m_i>>0$ and for the general curve $C$
cut out by $m_1 \; H_1,...,m_{n-1}H_{m-1},$ the restriction ${\cal E}|_C$ is
ample (nef).

\bigskip

Miyaoka [Mi 87] has shown that $X$ is not uniruled iff $\Omega^1_X$ is
generically $(H_1,...,H_{n-1})$-\underline{nef} for some $H_i$.

We obtain as a special case of (2.5):

\bigskip
{\bf Corollary} \hspace{0,3cm} {\it If $T_X$ is generically
$(H_1,...,H_{n-1})$-ample, then $X$ is rationally
generated.}

\bigskip

The converse of the corollary can be expected
but should be hard (at least with the $H_i$ big instead of ample).
It is much easier to look only at positivity of large
families of curves, not necessarily cut out by divisors. Also a relation to
stability is pointed out.

\bigskip

In practice the positivity of $T_X$ is difficult to handle. Much easier is
the anticanonical bundle $ -K_X = det \;  T_X.$  Therefore we shall investigate
positivity properties of $ -K_X$ in sect. 3 for varieties of negative
Kodaira dimension.

\section*{0. Notations around rational curves.}

Herewe collect for the convenience of the reader various concepts dealing
with rational curves. General references are [Ca 92], [KoMiMo 92].

\bigskip

{\bf(0.1)} A variety $X$ is called rationally connected if and only if two
general points on $X$ can be joined by a chain of rational curves.

\bigskip

{\bf(0.2)} Let $(C_t)$ be a covering family of rational curves on $X$.
Introduce the following equivalence relation: $ x \sim  y \Longleftrightarrow
x,y $ can be joined by a chain of $C_t's.$ Then there is a variety $Y$ and a
rational map $f \colon X \rightharpoonup Y$ such that for general
$x,y$ we have $f(x) = f(y) \Longleftrightarrow x  \sim y.$ The map $ f \colon X
\rightharpoonup Y $ is called the rational quotient of $X$ w.r.t. $(C_t)$
(see [Ca94] for
a proof, in a more general situation).

\bigskip

{\bf(0.3)} One can do the construction in (0.2) for all $(C_t)$ at the same time
and obtains the notion of a ''rational quotient'' of $X$.

More precisely define $x \approx y \Longleftrightarrow x,y $ can be joined
by a chain of rational curves.

Then there exists a variety $Z$ and a rational map $ g \colon X \rightharpoonup
Z $ such that for very general $x,y$ (i.e. outside a countable union of
subvarieties of $X$) we have $g(x) = g(y) \Longleftrightarrow x \approx y.$

An important property of the rational  maps $f \colon \rightharpoonup Y$
 and $g \colon \rightharpoonup Z$ is the almost holomorphicity, i.e. they
are holomorphic
and {\it proper} on a suitable Zariski dense open subset of $X.$

\bigskip

{\bf(0.4)} Let $X$ be a variety. Define inductively: $ X_o = X, \; X_i = $
rational quotient of $X_{i-1}.$

Then $X$ is said to be rationally generated if there exists $m$ such that
$dim \, X_m = 0.$ It is conjectured (but known only up to dimension 3,
[KoMiMo92]) that $X$ is
rationally connected iff $X$ is rationally generated.

Rational connectedness and rational generatedness are equivalent provided
the following holds:
\smallskip for every holomorphic surjective map $f \colon X \EP _1 $ with
$X$ a projective
manifold and $f$ having rationally connected fibers, there exists a
rational curve $C \subset
X$ with $f(C) = \EP _1.$ See [Ca95].

\section*{1. Ample subsheaves of tangent bundles}

In this section we shall prove

\bigskip

{\bf 1.1 Theorem} \hspace{0,3 cm}{\it Let $X$ be a projective manifold of
dimension $n$. Let $ {\cal E}  \subset {\cal T}_X $ be a locally free sheaf
of rank
$r$. If  ${\cal E}$ is ample and $ rk \; {\cal E} \ge n-2, $ then
$ X \simeq \EP_n $ and $ {\cal E} = {\cal T}_X $ or $ {\cal O}(1)^{\oplus r}$.}

\bigskip

Actually (1.1) should hold without any assumption on the rank. If
$ {\cal E} = {\cal T}_X,$ then (1.1) is nothing than Mori's theorem [Mo 79].
Also the other extremal case is known: if  $ rk \; {\cal E} = 1, $
then $ X \simeq \EP_n $ by J. Wahl [Wa 83].

\bigskip

The proof of (1.1) proceeds by induction on $ dim \; X $.  In the
induction step we will need the following

\bigskip

{\bf 1.2 Lemma} \hspace{0,3 cm}{\it Let $X$ be a n-dimensional projective
manifold,
$ \varphi \colon X \to Y $ a $\EP_k$-bundle $(k < n) $ of the form $X = \EP(V)$
with a vector bundle $V$ of rank $k+1$. Then the relative tangent sheaf
$ {\cal T}_{X|Y} $ does not contain an ample locally free subsheaf.}

\bigskip

{\bf Remark} \hspace{0,3 cm} Lemma 1.2 should only be a special case of a
much curve general fact: if $\varphi \colon X \to Y $ is a fiber space
(sufficiently smooth?), then ${\cal T}_{X|Y} = (\Omega^{1}_{X|Y})^*$ does not
contain an ample subsheaf unless $Y$ is a point and $X$ projective space.

\bigskip

{\bf Proof.} \hspace{0,3 cm} Suppose $ {\cal E} \subset {\cal T}_{X|Y} $ is
an ample
locally free subsheaf. Let $F$ be a fiber of $\varphi $. Then
$ {\cal E} |_F \subset {\cal T}_{X|Y} | F = {\cal T}_F.$
Hence either $ {\cal E} = {\cal T}_{X|Y} $ or ${\cal E} = {\cal
O}(1)^{\oplus r}_{|F}$
for all $F$.
In the first case, apply the Euler sequence

\bigskip

(E) $ \qquad 0 \to {\cal O} \to \varphi^*(V^*) \otimes {\cal O}_{\EP(V)}(1)
       \to {\cal T}_{X|Y} \to 0 $

\bigskip

to derive the ampleness of $\varphi^*(V^*) \otimes {\cal O}_{\EP(V)}(1) $
((E) does not split !) which is absurd (or argue as in the second case).
So let $ {\cal E}|_F = {\cal O}(1)|^{\oplus r}_F $ for all $F$.
Then

\bigskip

(*) $ \qquad {\cal E} \simeq \varphi^* ({\cal E}') \otimes {\cal O}_{\EP(V)}(1)$

\bigskip

with a bundle  $ {\cal E}'$ of rank $ r \; $ on $ \; Y $. Consider
again (E). We claim that the inclusion $ {\cal E} \hookrightarrow
{\cal T}_{X|Y} $
lifts to an inclusion $ {\cal E} \hookrightarrow \varphi^*(V^*) \otimes
{\cal O}_{\EP (V)} (1).$

=46or this we need to show that the canonical map
$$
   H^o({\cal E}^* \otimes \varphi^* (V^*) \otimes {\cal O}_{\EP(V)}(1))
   \to H^o ({\cal E}^* \otimes {\cal T}_{X|Y})
$$

is onto, hence that $H^1({\cal E}^*) = 0.$

By Leray's spectral sequence,
$$
    H^1 (X,{\cal E}^*) \simeq H^1(Y, {\cal E}'{}^* \otimes
    \varphi_*({\cal O}_{\EP(V)}(-1)) = 0.
$$

So we have $ {\cal E} = \varphi^*({\cal E}') \otimes {\cal O}_{\EP(V)}(1)
           \subset \varphi^*(V^*) \otimes {\cal O}_{\EP(V)}(1), $
therefore $ {\cal E}' \subset V^*.$  Consider the dual map
$$
    \alpha \colon V \to {\cal E}'{}^*,
$$

which is generically onto. Let $ {\cal S} = Im \; \alpha. $ Then

$    i : \EP({\cal S})   \hookrightarrow \EP(V) $ and

$$
     i^*({\cal E}) = i^*({\cal O}_{\EP(V)}(1) \otimes \varphi^*({\cal E}')) =
     {\cal O}_{\EP({\cal S})}(1) \otimes \varphi '{}^*({\cal E}'),
$$

  where  $ \varphi ' \colon \EP({\cal S}) \to Y $ is the projection.
  By construction, $ {\cal O}_{\EP({\cal S})}(1) \otimes \varphi'{}^*
  ({\cal E}) $ is ample, hence, taking det,

  $$
    {\cal O}_{\EP({\cal S})} (r) \otimes \varphi '{}^*(det {\cal E }') =
    {\cal O}_{\EP({\cal S} \otimes{{det \; {\cal E}'}\over {r}})} (r)
  $$

 is ample. Hence $ {\cal S} \otimes {{det \, {\cal E}'}\over {r}} $  is
 ample. Now take a general curve $C \subset Y $. Then $ {\cal E}'{}^* \otimes
 {{det {\cal E}'}\over {r}} | C $ is ample. But
 $ c_1({\cal E}'{}^* \otimes {{det {\cal E}'}\over {r}}) = 0, $ contradiction.

 \bigskip

 Now we come to the proof of (1.1). We treat the cases $ r=n, n-1$ and
 $n-2$ separately.

 \bigskip

{\bf(1.3)} Case $r = n. $

 The inclusion $ {\cal E} \subset {\cal T}_X $ yields $ det \, {\cal E} \subset
 - K_X$ . If $ det \, {\cal E} = - K_X, $ then $ {\cal E} = {\cal T}_X $ and
 the claim is Mori's theorem. So assume $ det \, {\cal E} \not=
 - K_X, $ hence $ K_X + det \, {\cal E} \not= {\cal O}_X,$ but $K_X + det \,
 {\cal E} \subset {\cal O}_X. $

 So $K_X + det \, {\cal E} $ is not nef. Now vector bundle adjunction theory
 [Fu 90, YZ 90] gives $X \simeq \EP_n $ and $ {\cal E} = {\cal O}(1)^n. $

 \bigskip

 {\bf(1.4)} Case $ r = n - 1. $

 Since all arguments in this case appear also in the more difficult case
 $ r =n-2, $ we just give a very rough sketch. $ {\cal E} \subset {\cal T}_X$
 yields $det \, {\cal E} \subset \land^{n-1}{\cal T}_X = \Omega^1_X \otimes -
 K_X, $ hence $ L = K_X + det \, {\cal E} \subset \Omega^1_X. $

 If $L$ is not nef, then all such pairs $(X, {\cal E})$  can be classified
  ([Fu 90, YZ 90]) and examined to prove our claim. So let $L$ be nef. Hence
 $$
      H^0(X, \Omega^1_X \otimes L^*) \not= 0
 $$

  for some nef line bundle $L$.

\bigskip

Now $H^o (\Omega^1_X \otimes L^*) \not= 0 $ implies that $X$ is not rationally
connected. However it is uniruled by [Mi 87], so one can consider the rational
quotient $ f \colon X \rightharpoonup Z $ w.r.t. some covering family of
rational
curves. Then the general compact fiber of $f$ is a projective space which
leads easily to a contradiction to the nefness of $L$.

\bigskip

{\bf (1.5)} Case $r = n-2. $

Taking $det$ of $ {\cal E} \hookrightarrow {\cal T}_X $ gives
$ det \, {\cal E} \hookrightarrow \land^{n-2}{\cal T}_X = \Omega^2_X \otimes
  - K_X $ hence

$$
     L := K_X + det \, {\cal E} \subset \Omega^2_X .
$$

\bigskip

{\bf (1.5.1)} We first assume that $L$ is not nef.

Then by [AM 95], $(X, {\cal E}) $ is one of the following.
\begin{itemize}
\item[(a)] There is a blow-up $ \varphi \colon X \to W $ of a finite set
           $ B \subset W $ of smooth points and an ample bundle ${\cal E}' $
           on $W$ with $ {\cal E} = \varphi^*({\cal E}') \otimes {\cal O}_X
           (-A), \qquad A = \varphi^{-1}(B) $
\item[(b)] $ X = \EP_n, {\cal E} $ splits
\item[(c)] $ X = Q_n $ and $ {\cal E} $ splits or is a twist of the
           spinor bundle
\item[(d)] $ X = \EP_2 \times \EP_2, \; {\cal E} = {\cal O}(1,1) \oplus
             {\cal O}(1,1)$
\item[(e)] $X$ is a del Pezzo manifold with $b_2 = 1, {\cal E} = {\cal O}(1)
           \oplus n-2 $
\item[(f)] $ X$ is a $ \EP_{n-1} $ - bundle or a quadric bundle over a
           curve $Y$
\item[(g)] $X$ is a $\EP_{n-2} $-bundle over a surface.
\end{itemize}

We comment on those cases separately.
\begin{itemize}
\item[(a)] We have $ {\cal E} = \varphi^*({\cal E}') \otimes {\cal O}_X
           (-A) \subset {\cal T}_X \subset \varphi^*({\cal T}_W), $
           so that $ {\cal E}' \otimes I_B \subset \varphi_*({\cal T}_X)
           \subset {\cal T}_W, $ hence
           $$
              {\cal E}' \subset {\cal T}_W .
           $$
           By induction on $ b_2(X) $ we may assume $W = \EP_n $ and
           $ {\cal E}' = {\cal O}(1)^{\otimes(n-2)}. $ But then
           $ {\cal E} = \varphi^*({\cal O}(1))^{n-2} \otimes {\cal O}_X(-A)$
           is clearly not ample, at most nef.
\item[(b)] Obviously $ {\cal E} = {\cal O}(1)^{\oplus(n-2)}.$
\item[(c)] If $ {\cal E} $ splits, then $T_{Q_n} $ would have a section
           vanishing on ample divisor which is not possible. So
           ${\cal E} $ is the twist $F(2)$ of a spinor bundle, see e.g.
           [PSW 90]. Therefore ${\cal E}$ has a section vanishing on
           a hyperplane section, so does $T_{Q_n} $  which would imply
           $ H^o(T_{Q_n}(-1)) \not= 0, $ contradiction.
\item[(d)] Again $ \EP_2 \times \EP_2 $ has no vector field vanishing on an
           ample divisor.
\item[(e)] Here $ K_X + det \; {\cal E} = {\cal O}_X(-1) $  and
           $ {\cal E} = {\cal O}_X(1)^{\oplus(n-2)}. $
           Now either apply [Wa 83], [MS 78] or apply the Fujita
           classification of del Pezzo manifolds or argue directly as
           follows. We have $H^o(X, {\cal T}_X(-1)) \not= 0. $ Choose
           $H \in |{\cal O}_X(1)| $ smooth and consider the sequence
           $$
             0 \to {\cal T}_H \to {\cal T}_X | H \to N_{H|X} = {\cal O}(1)
             \to 0 .
            $$
            If the induced map $ {\cal E}|_H \to N_H $ is non-zero, we obtain
            a map ${\cal O}(1)^{n-3} \hookrightarrow {\cal T}_H $,
            otherwise $ {\cal O}(1)^{n-2} \subset {\cal T}_H. $
            In both cases we get $ H \simeq \EP_{n-1} $ by induction, hence
            $X \simeq \EP_n.$
\item[(f)] Let $F$ be the general fiber of $ \varphi \colon X \to Y. $
           Then $ {\cal E}|_F \hookrightarrow {\cal T}_X |F $ and since the
           composition
           $$
             {\cal E}|_F \to {\cal N}_{F|X} = {\cal O}_F^a
            $$
            is zero, we obtain $ {\cal E}|_F \hookrightarrow {\cal T}_F $.
            Inductively, $ F \simeq \EP_{n-1}, $ hence $f$ cannot be a quadric
            bundle therefore must be a $\EP_{n-1}$-bundle; moreover
            $ {\cal E} \subset {\cal T}_{X|Y}. $
            But this contradicts (1.2).
\item[(g)] Finally $ \varphi $ is a $ \EP_{n-2}$ bundle.  As in (f),
           $ {\cal E} \subset {\cal T}_{X|Y},$ contradicting (1.2).
\end{itemize}

This finishes the case that $L$ is not nef.

\bigskip

{\bf (1.5.2)} $ L = K_X + det \, {\cal E} $ is nef.

So (*) $ \qquad H^o(X, \Omega^2_X \otimes L^*) \not= 0 $

with $L$ nef. This implies that $X$ is not rationally connected:
otherwise we find by [KoMiMo 92]
a family of rational curves $C_t \subset X $ with
$ T_X|C_t $ ample for $t$ general contradicting (*).

Hence we have a rational quotient $f \colon X \rightharpoonup Z,$ [Ca92]
and $ dim \, Z > 0 $ since $X$ is not rationally connected. Moreover $ dim \,
Z < dim \, X, $ since $X$ is uniruled by [Mi 87].
Since $f$ is
almost holomorphic, i.e. proper on an open set, it has some compact fibers.
The general fiber $F$ therefore must be $ \EP_k$ (since $ {\cal E}|_F
\hookrightarrow {\cal T}_F) $ with $k \ge n-2.$ Therefore $L_F$ is never nef,
contradiction.
\qed

\bigskip

{\bf(1.6)} We already mentioned Wahl's theorem several times: if $X$ is a
projective
manifold, $L$ an ample line bundle, $H^o(X,T_X \otimes L^*) \not= 0, $ then
$ X \simeq \EP_n. $ Previously Mori and Sumihiro [MS 78] proved this under
the additional assumption that $L$ is effective.

It should actually be possible to deduce Wahl's theorem from Mori-Sumihiro by
covering tricks in the following way. Take $ d \ge 2 $ with $ H^o(X, L^d)
\not= 0, $ and $ H^o(X,T_X \otimes L^*) \not= O.$ Suppose $ s \in H^o(X,L^d),
s \not= 0 $ with $ \{s=0\} $ smooth. Let $ f \colon Y \to X $ be the
associated cyclic cover so that $H^o(Y,f^*(L)) \not= 0. $ Then there is
an exact sequence

$$
  0 \to T_Y \to f^*T_X \to N_R \to 0
$$

where $N_R $ is supported on the ramification divisor $R \subset Y $ .
The map $ f^*(L) \hookrightarrow f^*(T_X) $ yields a map $ \gamma \colon
f^*(L) \to N_R. $ Now it can be shown that $ \gamma = 0 $, hence we
obtain $ f^*(L) \subset T_Y. $ By [MS 78], $Y \simeq \EP_n, $ moreover
$ f^*(L) = {\cal O}(1), $ so that  $deg \, f = 1 $ and $X \simeq \EP_n. $

Of course the zero set $ \{s = 0 \}$ could be non-smooth. If $ \{s = 0 \} $
is at least reduced, then the same arguments basically work, only that
$Y$ is merely normal. So one  should first prove [MS 78] in case $Y$ is normal.
The case where $ \{s = 0\} $ contains a non-reduced component seems more
complicated, one should proceed taking roots.

\section *{2.  The tangent bundle of rationally connected varieties}

{\bf 2.1 Definition} \hspace{0,3 cm} Let $X$ be a projective manifold. A
covering family
$(C_t)_{t \in T} \; (T $ compact) of curves is free if

(1) $T$ is irreducible

(2) $C_t$ is irreducible for $ t \in T $ general

(3) for all $A \subset X$ at least 2-codimensional there exists
    $ t \in T $ with $C_t \cap A = \emptyset $.

\bigskip We shall say that a family $(C_t)$ is $T_X$-ample, if $T_X \vert
C_t$ is
ample for general $t.$

\bigskip

{\bf 2.2 Theorem} \hspace{0,3 cm}{\it Let $X$ be a projective manifold,
$(C_t) $ a
free $T_X$-ample family. Then $X$ is rationally generated. Conversely, a
rationally
connected manifold has a free $T_X$-ample family of (rational) curves. }

\bigskip

{\bf Proof :} \hspace{0,3 cm} From the definition of a free family it is clear
that $X$ is uniruled, applying [MM 86]. Then we consider ''the'' rational
reduction $f \colon X \rightharpoonup Y $ and have $dim Y < dim X, \; Y $
smooth.
By construction, $Y$ is not uniruled.

We need to show that $dim Y = 0. $ So assume $ dim Y > 0 . $ By condition
(2) of (2.1), $f$ is holomorphic near $C_t$ for general $t$. Since
$ T_X|C_t $ is ample, we have $ dim \, f(C_t) > 0. $ Now consider the
generically surjective map ($t$ general)
$$
   \alpha \colon T_X|C_t \to f^*T_Y|C_t.
$$

Then $ Im \; \alpha $ is an ample sheaf, hence $ f^* T_Y|C_t$ is ample and
so does $T_Y|f(C_t). $ The family $C_t' = f(C_t) $ still covers
$Y$ with $T_Y|C_t' $ ample for general $t$, hence $ -K_Y.C_t' > 0 $ and
$Y$ is uniruled by [MM 86], contradiction.
The other direction follows from [KoMiMo92] using (2.5(2)).
\qed

\bigskip

Next we recall a notion due to Miyaoka.

\bigskip

{\bf 2.3 Definition} \hspace{0,3 cm} Let $X$ be a projective manifold,
$ n = dim \, X, H_1, \cdots , H_{n-1} $ ample divisors $X$. A vector bundle
$ E $ on $X$ is generically $ (H_1, \cdots , H_{n-1})$-ample (nef) iff for
$ m_1>>0, \cdots , m_{n-1} >> 0 $  and for general $C$ cut out by $m_1 \,
H_1, \cdots , m_{n-1} H_{n-1}$ the restriction $ E|C$ is ample (nef).

\bigskip

{\bf(2.4)} The same definition makes sense for $H_1, \cdots , H_{n-1}$
big, as long
as the general element cut out by $m_1 \; H_1, \cdots , m_{n-1}H_{n-1} $
is really a curve, i.e. if $ dim \bigcap\limits_{i=1}^{n-1} Bs(m_iH_i) \le
1.$

\bigskip

{\bf 2.5 Obvious Example} \hspace{0,3cm} (1) If $H_1, \cdots , H_{n-1} $
are very
ample, then the general $C_t$ cut out by $H_1, \cdots ,H_{n-1} $ deforms
to  build free families.

(2) If $C \subset X$ is a rational curve with $T_X \vert C$ ample, then $C$
moves in
a free $T_X$-ample family, provided $X$ is smooth near $C.$

\bigskip

{\bf(2.6)} Miyaoka has shown in [Mi 87] that $X$ is not uniruled iff $\Omega
^1 $ is generically $(H_1, \cdots , H_{n-1})$-nef for all ample
$ H_1, \cdots , H_{n-1} $ . We are now interested in the ''dual'' case.

\bigskip

{\bf 2.7 Theorem} \hspace{0,3 cm}{\it  If  $T_X$ is generically $(H_1, \cdots ,
H_{n-1})$-ample, then $X$ is rationally generated.}

\bigskip

This  follows from the more general theorem 2.2.

\bigskip

{\bf 2.8 Corollary} \hspace{0,3 cm} {\it Let $X$ be a smooth projective
3-fold, $(C_t)$
a free family of curves such that $T_X|C_t $ is ample for general $t$. Then
$X$ is rationally connected.}

\bigskip

{\bf Proof.} \hspace{0,3 cm} Rationally generated 3-fold are rationally
connected [KoMiMo 92].

\bigskip

{\bf 2.9 Remark} \hspace{0,3 cm}  One might expect that (2.2) remains true
if one assumes only the
ampleness of $T_X|C_t, t $ general, where $(C_t) $ is a family with the
following property: for $x,y \in X $ general, there is some $t$ such
$ x,y \in C_t $ (or joined by a chain of $C_t's$). However the proof
of (2.2) does not give this claim. Probably one should construct from such
a family a new, free family of curves.

\bigskip

{\bf 2.10 Problem} Let $X$ be rationally connected. Do there exist
$H_1, \cdots, H_{n-1} $ ample (or big), such that $T_X$ is generically
$ (H_1, \cdots , H_{n-1}) $ - ample? Is there at least a free family
$(C_t) $ such $T_X|C_t$ is ample, $t$ general?

\bigskip

{\bf (2.11)} We discuss (2.10) in several special cases.

\medskip

(1) If $n = dim \, X = 2 $, then (2.10) ''clearly'' holds.

\medskip

(2) Assume $T_X$ to be $(H_1, \cdots , H_{n-1})$ - semi-stable.

\medskip

Then $T_X|C $ is semi-stable for $C$ cut out by $ m_1 \, H_1, \cdots,
m_{n-1} \; H_{n-1} $ for  $m_i>>0 $ (Mehta-Ramanathan). By [Mi 87a], this
is equivalent to saying that
$ (T_X \otimes {{K_X}\over {n}})|C $ is nef (i.e. $S^nT_X \otimes K_X |C $
is nef) so if we know additionally that $-K_X \cdot H_1 \cdot ... \cdot
H_{n-1} > 0, $ (2.10) holds.

\medskip

So (2.10) has a positive answer if we can find $H_1, \cdots , H_{n-1} $
ample such that

\medskip

(a) $T_X $ is $(H_1, \cdots , H_{n-1}) $ semi-stable.

(b) $ H_1, \cdots , H_{n-1} $ is in the cone generated by classes of covering
    families of rational curves.

(since (b) implies $-K_X.H_1. \cdots . H_{n-1} > 0  \ !)$

\bigskip

In sect. 3 we shall discuss the relation of various cones of curves.

\medskip

(3) If $X$ is Fano, then (b) is always fulfilled.

It is expected that for $b_2 = 1 , T_X $ is always stable; a lot of
cases being checked in [PW 94], so that (2.10) should hold
for Fano manifolds with
$b_2 = 1 $ (and probably also for higher $b_2$).

\bigskip {\bf 2.12 Remarks} \hspace{0,3cm} (1) If $X$ is rationally
generated, does $X$
carry a free, $T_X$-ample family of curves? This could also be considered
as a converse of
the statement of (2.2).

(2) Is there a direct way to prove that a projective manifold carrying a free
$T_X$-ample family of rational curves is rationally connected?

(3) From (2) and (3) it would follow that rationally generated manifolds
are rationally connected. While (2) seems accessible, (3) is probably hard.

\section*{3. Various cones and the canonical bundle}

{\bf 3.1 Notation} \hspace{0,3 cm} $K_{nef}$ denotes the (closed) cone of nef
divisors, $K_{eff},$ that one effective divisors (always modulo numerical
equivalence, of course). We let $N_{nef}$ be the dual cone of $K_{eff}$
so that $ C \in N_1(X) $ is in $N_{nef}$ iff $D.C \ge 0 $ for all
$D \in K_{eff}$ or all effective divisors $D$ ) and call an element in
$N_{nef}$ a nef curve. The dual cone of $K_{nef}$
is of course $\overline{NE}(X)$, the closed cone of curves. Furthermore we let

\begin{itemize}
\item[a)] $ N_{rat} $ be the closed cone generated by classes $ [C_t] $,
           where $(C_t)$ is a covering family of rational curves, i.e. $C_t$ is
           a possibly singular irreducible rational curve for general $t$.
           Without requiring $C_t$ to be irreducible, but connected, with
all components
           of $C_t$ being rational, we call the resulting cone $\widetilde{N}_
           {rat}.$
\item[b)]  $N_{cov}$ be the analogous cone, omitting ''rational''
\item[c)]  $ N_{ci}$ be the closed cone generated by curves $H_1, \cdots ,
H_{n-1};
            H_i $ very ample ("ci" stands for complete intersection)
\item[d)] $N_{sc}$ be the closed cone generated by curves $H^{n-1}; \; H $
very ample  ("sc"
stands for special complete intersection).
\end{itemize}

\bigskip We have the obvious inclusions $N_{rat} \subset \tilde N_{rat}
\subset N_{cov}$
and $N_{sc} \subset N_{ci} \subset N_{cov}.$
\bigskip

{\bf 3.2 Problems} \hspace{0,3 cm} (cp. [DPS 95].

\begin{itemize}
\item[(1)] Is $ N_{nef} = N_{cov}? $ (Clearly $N_{cov} \subset N_{nef}).$
\item[(2)] Assume $X$ Fano. Is $N_{rat} = N_{nef}$ or at least
           $N_{rat} = N_{cov} $?
\item[(3)] Assume $X$ rationally connected or Fano. Is
           $N_{rat} \cap N_{ci} \not= \{0\} $?
\item[(4)] Assume $X$ rationally connected or Fano. Is
           $N_{rat} \not= \widetilde{N}_{rat}$?
\end{itemize}

\bigskip

{\bf 3.3 Remarks}

\begin{itemize}
\item[(1)] This problem was posed in [DPS 95]. It "essentially" (up to a
limit process) means that given
           an nef curve $C \subset X, $ i.e. $D.C \ge 0 $ for all effective
           $D$, then there should exist a covering family $(C'_t)$ such
           that $C'_t \equiv \alpha \; C_t $ for some $ \alpha \in \EQ_t $.
           In particular, if every effective divisor is nef, e.g.  $T_X$
           is nef, then every curve $C$ should be ''numerically movable''
           up to a multiple. For surfaces this is clear. This would give
           new evidence for the conjecture that Fano manifolds with
           $T_X$ nef are homogeneous ([CP 91, DPS 94]).

\item[(2)] Assume $X$ Fano. The intention is to investigate ''how much of
the geometry
is dictated by the moving rational curves.'' The cone theorem says that
$\overline{NE}(X)$ is already generated by the classes of rational curves.
So it seems very natural to suspect that the cone of ''covering curves''
is generated by the cone of ''covering rational curves''. As an application,
a divisor $D$ would be effective iff $D.C \ge 0 $ for all rational curves
$C$ with $T_X | C$ nef. If $T_X$ itself is nef, this follows from
[CP 91] . Note that for any $X$ the equality
$N_{rat} = N_{nef}$ implies in particular that
$ - K_X \in K_{eff}$. So (2) cannot hold for every rationally connected
manifold. But it could hold for every rationally connected manifold with
$ -K_X \in K_{eff}$.

\item[(3)] This is motivated by stability, see (2.10), (2.11).
\end{itemize}

\bigskip

{\bf 3.4 Proposition} \hspace{0,3 cm} {\it Let $X$ be a del Pezzo surface. Then
$N_{rat} = N_{nef}.$}

\bigskip

{\bf Proof}. \hspace{0,3 cm} Since $dim \, X = 2, $ we have $N_{amp} =
K_{amp}. $
Now $K_{amp}$ is generated by $ \varphi^*H', $ where $ \varphi \colon
X \to X' $ is the contraction of a codimension 1-face $R$ in
$ \overline{NE}(X) $ and $H'$ is ample on $X'$. Moreover, $X$ being del
Pezzo, $X' = \EP_1 $ or $X' = \EP_2 $ .

If $X' = \EP_1, $ then $R$ is given by $\varphi^{-1}$ (point $p$), if
$X' = \EP_2$, then $R$ is given by $\varphi^{-1}$ (general line $\ell$).
Both $ \varphi^{-1}(p)$ and $\varphi^{-1}(\ell) $ are nef curves, finishing
the proof.

\bigskip

Concerning 3.2.(3) we prove:

\bigskip

{\bf 3.5 Theorem} \hspace{0,3 cm}{\it  Let $X$ be a Fano manifold.}

Then $ N_{sc} \subset \widetilde{N}_{rat}. $

\bigskip

{\bf Proof.} \hspace{0,3 cm} By the cone theorem (as in (3.4)) it is sufficient
to show the following

\bigskip

Claim (*) \hspace{0,3 cm} Let $ \varphi \colon X \to Y$
be the contraction ofa codimension 1 face $R$ and $H$ ample on $Y$. Then
$ \varphi^*(H^{n-1}) \in \widetilde{N}_{rat}. $

a) If $ dim \, Y < dim \, X, $ we proceed as follows. First note
if $dim \, Y \le n-2, $ then $H^{n-1} = 0 $, so the claim is obvious.
So let $dim \, Y = n-1. $ Then $H^{n-1}$ can be thought of as a  general
point of $Y$ (up to a positive rational multiple of $H^{n-1}$), so that
$ \varphi^*(H^{n-1}) \equiv $ general fiber of $ \varphi $, which is a smooth
moving $ \EP_1, $ hence in $N_{rat}$, proving (*) in case $ dim \, Y < dim
\, X$ .

b) $ \varphi $ is birational.

Note that $ \rho(Y) = 1 $ and that $Y$ is rationally connected (but may
be singular).

=46irst let us prove the following assertion

\bigskip (**) Let $C \subset Y $ be an irreducible curve not contained in
the degeneracy
set in $Y$ of $ \varphi \colon X \to Y. $

Then there is a - up to numerical equivalence -
unique effective curve $ \widetilde{C} \subset X $ such that

(1) $ \varphi_*(\widetilde{C}) = C $

(2) $ \widetilde{C}.E_i = 0 $ for all exceptional components of $\varphi$.

\bigskip

In fact, let $ \widehat{C} \subset X $ be the strict transform of $C$ in $X$.
We make the following ansatz
$$
   \widetilde{C} = \widehat{C} + \sum \lambda_i C_i ,
$$

 where $C_i$ are the extremal rational curves contracted by  $\varphi$. We
 choose the $ \lambda_i $ in such a way that $ \widetilde{C}.E_j = 0 $
 for all divisors $E_j$ contracted by $ \varphi $. Then the $ \lambda_i$
 are uniquely determined. What remains to be shown is $ \lambda_i \ge 0.$
 Take a basis of the $\varphi$-ample divisors $D$, with $\varphi^*(C).D = 0$,
 say $D_1, \cdots , D_{\rho-1}.$ Every $D_j$ is of the form $D_j = - \sum
 \kappa_\mu E_\mu $ with $\kappa_\mu \ge 0 $ by Lemma 3.7. Hence we have

\bigskip

$ D_j.\sum \lambda_i C_i = D_j.\widetilde{C} - D_j.\widehat{C} -D_j.\widehat{C} \sum \kappa_k E_k.\widehat{C} \geq 0.$

 \bigskip

 Therefore $ \sum \lambda_iC_i
   \in \overline{NE}(X|Y),$ i.e. all $\lambda_i \ge 0. $ So (**) holds.

\bigskip

Now the verification of (*) in the birational case is easy; choose a covering
family $(C_t)$ of rational curves in $Y$. Then, applying (**), we have for $t$
general an effective curve $\widetilde{C}_t \subset X$ with
$\varphi_*(\widetilde
{C}_t) = C_t$ and $\widetilde{C}_t.E_i = 0 $ for all $\varphi-$exceptional
divisors $E_i \subset X $. We can choose $ \widetilde{C}_t $ rational because
$\varphi$ is the contraction of an extremal face.
Now let $H$ be the ample generator of $ Pic(Y) \simeq \EZ.$
Take $\mu \in \EN$ such that $\mu H $ is very ample so that we find a smooth
curve $\Gamma \equiv (\mu H)^{n-1}. $ Then $\Gamma \equiv mC_t$ for some
$ m \in \EQ_+$ since $\rho(Y) = 1. $ It follows that $\Gamma \equiv m \widetilde
{C}_t, $ where $\widetilde{\Gamma}$ is again constructed by (**). On the other
hand, $\widetilde{\Gamma} \equiv \varphi^*(H)^{n-1}.$ Hence the
$(\widetilde{C}_t)$ yield (taking closure) a covering family of rational
curves, numerically of the form $ \varphi^*(H^{n-1}),$ so that
$ \varphi^*(H^{n-1}) \in \widetilde{N}_{rat} (X)$.

There are two problems left: first, our $\widetilde C$ might not be
connected and second,
the $\lambda_i$ might only be rational and anyway no multiplicities are
allowed for the
curves in $N_{sc}.$
To get around with these difficulties we first choose $C$ to be a general
member of a large family
of rational curves in $Y,$ say a free family. Then choose a positive
integer $k$ such that
$k \lambda_i$ are integers for all $i.$ Now substitute $kC$ by general
deformations $C_j$ of $C.$
The $\lambda_j$ will then be the same for $C$ and $C_j.$ Since all fibers
of extremal contractions
are uniruled by [Ka92], we can substitute the $\lambda_i \widetilde C_i$ by
rational curves
$C_{i,1}, \ldots, C_{i,k\lambda_i} $ homologous to multiples of $\widetilde
C_i$ such that
the new $\widetilde C = \bigcup \hat C_j \cup \bigcup C_{i,j}$ is connected.

\bigskip

{\bf 3.6 Remark.} \hspace{0,3 cm} In order to prove $N_{sc} \subset
N_{rat}$ in (3.5) we would need the following.

Let $Y$ be a singular Fano variety (terminal singularities), $ \rho(Y) = 1 ,
A \subset Y $ an algebraic subvariety of codimension $\ge 2,$ then
there exists a covering family $(C_t)$ of rational curves, $C_t$ irreducible
for $t$ general, such that $ A \cap C_t = \emptyset $ for $t$ general.

\bigskip

{\bf 3.7 Lemma} \hspace{0,3 cm} {\it Let $X$ be a projective manifold,
$f \colon X \to Y $ be a birational morphism to a normal projective
variety $Y$. Let $E_1, \ldots ,E_k $ be the codimension 1 components of the
exceptional set of $f$. Let $D$ be $f$-ample. Then

$$
    D = - \sum \kappa _{\mu} \, E_{\mu} + f^*(L)
$$

with $\kappa_{\mu} \ge 0$ for some line bundle $L$ on $Y.$}

\bigskip

{\bf Proof.} \hspace{0,3 cm} Of course an equation as above always exists;
we only need to show $\kappa_{\mu} \ge 0.$ Taking hyperplane sections we
reduce to the case $ dim \; X = 2.$  Let $ B = \{ \mu | \kappa_{\mu} \le 0
\}$ and $A$ be its
complement; assume $B \not= \emptyset.$ Since $D$ is $f$-ample, we have

$$
   \sum\limits_{\mu} \, \kappa_{\mu} \, E_{\mu} \cdot \sum\limits_
   {\nu \in B} (-\kappa_{\nu}) E_{\nu} < 0.
$$

Hence
$$
    \sum\limits_{\mu \in A\atop \nu \in B}
        \kappa_{\mu}(-\kappa_{\nu}) E_{\mu} \cdot E_{\nu} <
        \sum\limits_{\mu, \nu \in B} \kappa_{\mu} \kappa_{\nu}
        E_{\mu} \cdot E_{\nu} .
$$

        Now the right hand side of the last equation is negative, the
intersection matrix
        $ (E_{\mu} \cdot E_{\nu}) $ being negative definite, while the left
hand side is
        obviously non-negative, contradiction.
\qed

\bigskip

We now relate the canonical bundle to various cones. The question we have
mind is the following: how negative is the canonical bundle $K_X$ of a manifold
$X$ with $\kappa(X) = - \infty? $ In dimension $> 3 $ it is even unknown
whether $K_X$ is not nef, so we will mainly restrict ourselves to $dim \, X=3$.

We introduce the following notation:

$$
    \overline{NE}_-(X) = \{ C \in \overline{NE}(X) | K_X \cdot C < 0 \}.
$$

Then vaguely speaking, we ask for the size of $\overline{NE}_-(X). $ E.g.:

 \noindent

{\bf 3.8  Problem} \hspace{0,3 cm} Does $ \kappa(X) = - \infty $
imply that $\overline{NE}_-(X) \cap N_{clt} \not= \{0\}? $

In other words: are there $H_1, \cdots , H_{n-1} $ ample such that
$ K_X \cdot H_1 . \cdots . H_{n-1} < 0$? 

\bigskip

A weaker problem would ask for $H_i$ only big. Note that this problem is
related to the stability considerations in \PS 2 .
\bigskip We shall now consider the case $dim X = 3 . $

\bigskip

{\bf 3.9 Proposition} \hspace{0,3 cm} {\it
Let $X$ be a normal ($\EQ$-factorial)
threefold with at most terminal singularities, $\varphi \colon X \to Y $
the contraction of an extremal ray with $ dim \, Y \le 2 $. Then there exists an
ample divisor $H$ on $X$ with $ K_X.H^2 < 0 . $ }

\bigskip

{\bf Proof.} \hspace{0,3 cm} If $ \rho(X) = 1, $ the claim is obvious, so
let $\rho(X) \ge 2, $ i.e. dim $Y>0$. Fix $H_o \in Pic(X) $ ample. Choose
any ample divisor $L$ on $Y$ and let

$$
  H = H_o + \varphi^*(mL),
$$

$m>>0$ to be specified in a moment. $H$ is ample and

$$
   K_X \cdot H^2 = K_X \cdot H_o^2 + 2 K_X \cdot H_o \cdot \varphi^*(mL) +
   K_X \cdot \varphi^*(mL)^2 .
$$

If $ dim \, Y=2, $ then $K_X \cdot \varphi^*(L)^2 < 0 $ since $- K_X$ is
$\varphi$-
ample, hence $K_X \cdot H^2 < 0 $ for $m>>0.$

If $ dim \, Y = 1 $, then $L^2 = 0 $, but now $K_X \cdot H_o \cdot
   \varphi^* (mL) < 0, $ hence again the claim for $m>>0. $

\bigskip

Note that we only used that $-K_X $ is ample on the general fiber instead
of $\varphi$ being an extremal contraction.

\bigskip

{\bf 3.10 Proposition} \hspace{0,3 cm} {\it Let $X,Y $ be normal
projective ($ \EQ$-factorial) terminal 3-folds, $\varphi \colon X \to Y $
a divisorial extremal contraction. Assume there exists $H'$ ample on
$Y$ with $K_Y \cdot H'{}^2 < 0. $ Then there is $H$ ample on $X$ with
$K_X \cdot H^2 < 0. $}

\bigskip

{\bf Proof.} \hspace{0,3 cm} Choose $H_o$ ample on $X$ and let
$H = H_o + \varphi^*(mH).$ Then $H$ is ample and
$K_X \cdot H_o^2 + 2 K_X \cdot H_o \cdot \varphi^*(mH) + K_X \cdot
\varphi^*(mH)^2. $
Since $K_X \cdot \varphi^*(H)^2 = K_Y \cdot H^2 < 0, $ we have
$K_X \cdot H^2 < 0 $ for $ m >> 0 . $

\bigskip

{\bf 3.11 Proposition} \hspace{0,3 cm}{\it Let $X$ be as in (3.10). Let
$$
  \begin{array}{llllll}
   & X &- - -\to  &  X^+ &  \\
   & & & & & \\
  \varphi  & \hspace*{0,7 cm}\searrow  & \hspace*{1,5cm} \swarrow
&\hspace*{0,7 cm}\varphi^+ \\
  & & & & & \\
  & & \hspace*{0,5 cm}  Y & &
  \end{array}  \qquad \qquad \mbox{be a flip}
  $$

  so $-K_X $ is $\varphi$-ample and $K_X+ $ is $\varphi^+$-ample.

  Let $H^+$ be effective on $X^+$ whith $ H^{+2} \in NE(X).$ Let
  $H$ be its strict transform on $X$. If $ K_X \cdot H^2 < 0 $, then
  $K_{X^+} \cdot H^{+2} < 0. $}

  {\bf Proof.} \hspace{0,3 cm}  Write

  \bigskip

  $H^{+2} = \sum \lambda_i\, C_i^+  + \sum \mu_j \cdot A^+_j, \qquad
  \lambda_j > 0 , \quad \mu_j \ge 0,  $

\bigskip

  where the $A_j^+ $ are exactly the exceptional curves for $\varphi^+$ .

\bigskip

Analogously: $ H^2 = \sum \; \tau_i \, C_i + \sum \; \nu_j \, A_j . $

(it is clear that $H^2 \in NE(X)$ ! ). We arrange things, such that
$C_i$ is the strict tansform of $C_i^+$. Since
$K_{X^+}.A_j^+ > 0 $ and $K_X.A_j < 0 $ it suffices to show

(*) $ \quad K_X.C_i \ge K_{X^+} .C_i^+.$

Let $C = C_i, \; C^+ = C_i^+ $ for simplicity. Let

$$
   \begin{array}{llllll}
   & & \hspace*{0,7 cm}  \hat{X} & & &  \\
   & & & & & \\
   \pi & \swarrow & &  \searrow  \sigma \\
   & & & & & \\
   X &  &- - -\to & & X^+
   \end{array}
$$

\bigskip

be a common desingularisation. By [KMM 87, 5-1-11] we have

$$
  K_{\widehat{X}} = \sigma^*(K_{X^+}) + \sum \; a_i^+ \; E_i
$$
$$
  K_{\widehat{X}} = \pi^*(K_X) + \sum \; a_i \; E_i
$$

and $a_i \le a_i^+ $ .

Therefore (*) follows by intersecting with $ \widehat{C},$ the
strict transform of $C_i$ resp. $C^+ $ in $\widehat{X} $.

\bigskip

{\bf 3.12 Remarks} \hspace{0,3 cm} Of course we would like to show that
given a 3-fold
$X$ with $ \kappa (X) = - \infty$, then there exists $H$ ample (or at least big)
such that $K_X.H^2 < 0.$ The strategy would be to apply the Minimal Model
Program (MMP) to finally abtain a contraction $X' \to Y $ with $dim Y \le
2$, then
apply (3.9) and go back successively. However the flips cause trouble
because the inequalities work in the wrong direction. In fact, in the
situation of (3.11) we would need $H^{+2} \in NE(X) $ and $ (H^+)^2 \cap $
(the exceptional locus of the flip) $ = \emptyset $. If say $H^+$ is big
and globally generated, this could be achieved. However $H$ would in
general no longer be globally generated, so we cannot proceed.

In special situations however, the problem has a positive solution:

\begin{itemize}
\item[a)] there exists $\varphi \colon X \to Y, \; dim \, Y \le 2, -K_X $
          is $\varphi$-ample or at least $-K_X$ is
          ample on the general fiber,
\item[b)] the Albanese dimension of $X$ is 2 (as a consequence of (a)),
\item[c)] $ | - mK_X| \not= \emptyset $ for some $ m \in \EN $ (obvious)
\item[d)] there is a divisorial modification $\varphi \colon X \to X' $
          and $X'$ is of type a) (e.g. no flips occur in MMP).
\end{itemize}

\bigskip

It is clear from Prop. 3.9 that in general we have at least:

\bigskip

{\bf 3.13 Proposition} \hspace{0,3 cm} {\it Let $X$ be a normal projective
3-fold,
$\EQ$-factorial, terminal with $\kappa (X) = - \infty$. Then there exists a
covering family $(C_t)$ of curves, $C_t$ irreducible with \underline{ample}
normal bundle for $t$ general, such that $K_X.C_t < 0. $}

\bigskip

Uniruledness implies only the existence of a family $(C_t)$ (however rational)
whose normal bundles are nef ($t$ general).

\end{document}